\documentclass[11pt, a4paper]{article}
\usepackage{amsmath,amsthm,amssymb,amsfonts}
\usepackage{a4wide}
\newtheorem{theorem}{Theorem}
\newtheorem{lemma}[theorem]{Lemma}

\theoremstyle{remark}

\newtheorem{remark}[theorem]{Remark}

\def\Xint#1{\mathchoice
{\XXint\displaystyle\textstyle{#1}}%
{\XXint\textstyle\scriptstyle{#1}}%
{\XXint\scriptstyle\scriptscriptstyle{#1}}%
{\XXint\scriptscriptstyle\scriptscriptstyle{#1}}%
\!\int}
\def\XXint#1#2#3{{\setbox0=\hbox{$#1{#2#3}{\int}$ }
\vcenter{\hbox{$#2#3$ }}\kern-.6\wd0}}

\def\dashint{\Xint-}

\usepackage{graphicx}
\usepackage[usenames,dvipsnames]{color}
\usepackage[colorlinks=true, pdfstartview=FitV, linkcolor=black, citecolor=blue, urlcolor=blue]{hyperref}


\newcommand{\e}{{\mathbb E}}

\newcommand{\N}{{\mathbb N}}

\newcommand{\R}{{\mathbb R}}

\newcommand{\mO}{\mathcal{O}}

\newcommand{\mB}{\mathcal{B}}

\newcommand{\pa}{{\partial}}
\newcommand{\na}{{\nabla}}

\newcommand{\eps}{{\varepsilon}}

\def\div{\hbox{div \!}}

\def\mspace{\medskip \noindent}

\title{A simple justification of effective models \\ for conducting or fluid media with dilute spherical inclusions} 
\author{David G\'erard-Varet}

\begin{document}
\maketitle

\begin{abstract}
We present a gentle approach to the  justification of  effective media approximations, for PDE's set outside the union of  $n \gg 1$ spheres with low volume fraction. To illustrate our approach, we consider three classical examples: the derivation of the so-called {\em strange term}, made popular by Cioranescu and Murat, the derivation of the Brinkman term in the Stokes equation, and  a scalar analogue of the effective viscosity problem. Under some  separation assumption on the spheres, valid for periodic and random distributions of the centers, we recover  effective models as $n \rightarrow +\infty$ by simple arguments.   
 \end{abstract}

\section{Introduction}
Let $B_{i,n} = B(x_{i,n}, r_n)$, $1 \le i \le n$,  a collection of disjoint balls, included in a compact subset  $K$ of $\R^3$. We assume convergence of the empirical measure
\begin{equation} \label{empirical} \tag{H0} \frac{1}{n} \sum_{i=1}^n \delta_{x_{i,n}} \rightarrow \rho(x) dx, \quad n \rightarrow +\infty 
\end{equation}
where $\rho$ is a bounded density with support $K_\rho \subset K$. Moreover, we  assume that the volume fraction of the balls 
\begin{equation} \label{def_fraction}
 \lambda_n  := \frac{4\pi}{3} \frac{n r_n^3}{|K_\rho|}  
 \end{equation}
is small uniformly in $n$  (and in some cases vanishes as $n$ goes to infinity). 

\mspace
Our concern is the asymptotics of elliptic  PDE's of Laplace or Stokes type, in the domain $\Omega_n  := \R^3 \setminus \cup_i B_{i,n}$: 
\begin{equation} \label{general} 
L u_n = g  \quad  \text{in} \:  \Omega_n,
 \end{equation}
completed with  boundary conditions at each ball $B_i$ and a decay condition at infinity.  Such PDE's arise in various models, taken from various areas: electrostatics, optics, fluid mechanics, heat conduction and many more.  Our main three examples will be: 
\begin{description}
\item[i)] Diffusion in domains with holes: 
\begin{equation} \label{Diff}
 -\Delta u_n = g \quad  \text{in} \:  \Omega_n, \quad 
u_n\vert_{B_{i,n}} = 0, \quad 1 \le i \le n. 
\end{equation}
\item[ii)] Drag in fluid flows around obstacles: 
\begin{equation} \label{Sto}
 - \Delta u_n + \na p_n = g, \quad \div u_n = 0  \quad  \text{in} \:  \Omega_n, \quad
 u_n\vert_{B_{i,n}} = 0, \quad 1 \le i \le n. 
 \end{equation}
\item[iii)] Permittivity in composites with perfectly conducting inclusions : 
\begin{equation} \label{Con}
 - \Delta u_n  = g \quad  \text{in} \:  \Omega_n, \quad
 u_n\vert_{B_{i,n}} = u_{i,n}, \quad u_{i,n} \in \R, \quad \int_{\pa B_i} \pa_\nu u_n   d\sigma = - \int_{B_i} g,  \quad 1 \le i \le n, 
 \end{equation}
\end{description}
with $\nu$ the unit normal vector pointing outward. 

\mspace
Just like \eqref{Diff}  has  \eqref{Sto} for  fluid counterpart, \eqref{Con} is the scalar version of the effective viscosity problem 
\begin{equation} \label{Einstein}
\begin{aligned} 
& - \Delta u_n + \na p_n = g, \quad \div u_n = 0  \quad  \text{in} \:  \Omega_n,  \quad u_n\vert_{B_{i,n}} = u_{i,n}  + \omega_{i,n} \times (x - x_{i,n})  \quad u_{i,n}, \omega_{i,n} \in \R^3,  \\
&  \int_{\pa B_i} (2D(u_n)\nu  - p \nu) d\sigma = -\int_{B_i} g, \\
 & \int_{\pa B_i} (2D(u_n)\nu  - p \nu)(x) \times (x - x_i) d\sigma(x) = -\int_{B_i} g(x) \times (x-x_i) dx,  \quad 1 \le i \le n,  
 \end{aligned}
 \end{equation}
that could be included in our analysis as well. Other problems fit a similar framework, such as the propagation of waves in bubbly fluids \cite{MR3690650}. 

\mspace
Back to the general formulation \eqref{general}, the point is to determine if, for large $n$ and small volume fraction,  the solution $u_n$ is close to the solution $u$ of an {\em effective model}: 
\begin{equation} \label{general_effective} 
L u + L_{e} u  = 0 \quad \text{in} \: \R^3 
\end{equation}
where the extra {\em effective operator} $L_{e}$ reflects some average effect of the spheres. Moreover, the hope is to express $L_{e}$  in terms of  macroscopic characteristics of the spheres distribution, like the limit density $\rho$ and the limit volume fraction $\lambda := \lim_n \lambda_n$ (when non zero). This hope is supported by the following formal reasoning. If the  volume fraction of spherical inclusions is small, one can expect the average distance between the spheres to be large compared to their radius. Therefore,  the interaction of the spheres should be negligible, and their contribution to the effective operator should be additive. 
Furthermore, if  the volume fraction vanishes as $n$ grows, the spheres should be well approximated by points. By combining both arguments, the leading effect of the spheres  at large $n$ should be through  the empirical measure. In the limit $n \rightarrow +\infty$, one should then recover $L_e = L_e[\rho]$. 

\mspace
In this spirit, early formal calculations of effective models were done by Maxwell Garnett \cite{MaxwellGarnett}, Clausius \cite{Cla}, Mossotti \cite{Mos}, to describe heterogeneous electromagnetic media. Similar computation of the effective viscosity of dilute suspensions was done by Einstein \cite{Ein}. For the drag generated by obstacles in a fluid flow, see Brinkman \cite{Brinkman1949}. Transposed to  our three examples, these calculations lead to the following lose statements: 

\begin{description} 
\item[i)] For system \eqref{Diff},  the critical scaling is $r_n = \frac{1}{n}$, with convergence of $u_n$ to  the solution  $u$ of 
\begin{equation} \label{Diff_eff}
-\Delta u + 4 \pi  \rho u = g  \quad \text{in }  \R^3. 
\end{equation}
\item[ii)] For system \eqref{Sto},  the critical scaling is  $r_n = \frac{1}{n}$, with  convergence of $u_n$ to  the solution  $u$ of 
\begin{equation} \label{Sto_eff}  
-\Delta u+ \na p + 6 \pi \rho u = g, \quad \div u = 0 \text{ in } \R^3. 
\end{equation}
\item[iii)] For system \eqref{Con}, the critical scaling is when  $\lambda := \lim_n \lambda_n$ is small but non zero. One has in this case $u_n = u + o(\lambda)$ uniformly in $n$, with  
\begin{equation} \label{Con_eff}
- \div \big( (1+3\lambda|K_\rho| \rho)  \na u \big) = g  \quad \text{in $\R^3$. } 
 \end{equation}
\end{description}
By {\em critical scaling}, we mean that other scalings would lead to trivial limits: typically, in system \eqref{Diff}, one has $u = 0$ if $r_n \gg \frac{1}{n}$, while it solves the original Laplace equation if $r_n \ll \frac{1}{n}$.  

\mspace
Turning these assertions into rigorous mathematical statements has attracted a lot of attention since 1970. Pioneering results on  \eqref{Diff} are due to Hruslov and Marchenko \cite{MR0601222}, Cioranescu and Murat \cite{MR1493040}, Papanicolaou and Varadhan \cite{MR609184}, or Ozawa \cite{MR727440}. A more abstract approach through $\Gamma$-convergence was completed in \cite{MR648338,MR640775}. For the most up to date statements, we refer to the nice study by Giunti, H\"ofer and Vel\'azquez \cite{MR3915491}, and to the bibliography therein. 
On the asymptotics of \eqref{Sto}, one can mention the work of Allaire \cite{MR1079189}, as well as improvements by Desvillettes, Golse and Ricci\cite{MR2398959}, Hillairet \cite{MR3851058}, or H\"ofer \cite{MR3863071}. As regards the analysis of \eqref{Con} and of the more involved fluid problem \eqref{Einstein}, rigorous results were obtained  in Levy and Sanchez-Palencia \cite{MR813656,MR813657}, Ammari et al \cite{MR3025042}, Haines and Mazzucato  \cite{MR2982744}, as well as in the recent papers by Hillairet and Wu \cite{HiWu}, Niethammer and Schubert \cite{NiSc} , or the author and Hillairet \cite{DGV_MH}. These references are by no mean exhaustive: one could further cite \cite{MR2094523,MR2504035,MR3449757,Mec,MR3458165,Glor}  and many more on related problems. 

\mspace
As explained before, the derivation of the limit systems \eqref{Diff_eff}, \eqref{Sto_eff} and \eqref{Con_eff} is based on neglecting  the interaction between the spheres. It requires some separation assumption on the centers $x_{i,n}$ of the spheres.   The challenge  is to show convergence under the mildest possible separation assumption. This is what differentiates the works mentioned above. Three types of assumptions emerge from the literature: 
\begin{itemize}
\item the most stringent one is that the centers are periodically distributed over a mesh of typical length $n^{-1/3}$. This is assumed for instance in  \cite{MR1493040}, \cite{MR1079189} or in the studies  \cite{MR813656,MR813657} as well as in  \cite{MR2982744}. 
 Note that in such a case, the limit density $\rho$ is constant over its support: $\rho = |K_\rho|^{-1} 1_{K_\rho}$.   
\item Other studies relax the periodicity assumption, but keep in addition to \eqref{empirical} a condition on the minimal distance between the centers: 
\begin{equation} \label{minimal_distance}
 \inf_{i \neq j} |x_{i,n} - x_{j,n}| \ge c n^{-1/3}. 
 \end{equation}
This is the case in  \cite{MR727440}, or in \cite{MR2398959}, and as far as we know in all works on the effective permittivity problem \eqref{Con}  or on the effective viscosity problem \eqref{Einstein}. Although it is an improvement to the periodic case, such assumption is not fully satisfactory: in particular, it is not satisfied if the centers are drawn randomly and independently  with law $\rho(x) dx$. Indeed, in this latter case, it is known that the typical minimal distance scales like $n^{-2/3}$.  
\item Eventually, the asymptotics of system \eqref{Diff} or \eqref{Sto} may be established under weaker assumptions, which apply when the centers $x_{i,n}$ are given by i.i.d. random variables, or derived from reasonable stationary point processes. A big step in that direction was made in  \cite{MR609184}: roughly\footnote{actually,  \cite{MR609184} treats the time dependent version of the equations},  the  solution $u_n$ of \eqref{Diff} tends to the solution  $u$ of \eqref{Diff_eff} under the weak separation assumption 
\begin{equation} \label{weakseparation} 
\frac{1}{n^2} \sum_{i \neq j} \frac{1}{|x_{i,n} - x_{j,n}|} \rightarrow \int_{\R^3 \times \R^3} \frac{1}{|x-y|} \rho(x) \rho(y) dx dy, \quad n \rightarrow +\infty
\end{equation}
It is  shown in  \cite{MR609184}  that this convergence holds in probability in the i.i.d. case. Recently, almost sure convergence of $u_n$ to $u$ was proved in  \cite{MR3915491} for a large class of stationary point processes, allowing for random radii of the spheres. See also  \cite{GiuHof} for the analogue  in the Stokes case.    
\end{itemize}
The works mentioned above may rely on very different techniques. Still, key arguments are often difficult and based on explicit constructions. For instance, the work of Cioranescu and Murat, as well as  the refined analyses  \cite{MR1079189, MR3915491},  are  based on the construction of so-called  {\em correctors},  which may be quite technical notably in the Stokes case.  Another instance is \cite{MR609184}, where the proof relies on the explicit representation of $u_n$ with the Feynman-Kac formula, and on non-trivial manipulations of this formula. In particular, none of these methods really exploits the fact that the formal limit equation is already known. The goal of this paper is to present a softer and shorter approach to these asymptotic problems. This approach relies on separation hypotheses in the same spirit  as \eqref{weakseparation}, although slightly stronger. In the context of problem \eqref{Diff}, resp.  \eqref{Sto}, our assumptions \eqref{H1}-\eqref{H2}, resp.  \eqref{H1}-\eqref{H2'}, see next section,  cover the periodic and random i.i.d. cases. In the context of \eqref{Con}, we derive the effective permittivity formula under a couple of conditions \eqref{A1}-\eqref{A2} that is more general than \eqref{minimal_distance}, and  to our knowledge is new. Still, we insist that the point of this paper is not the novelty of our statements, for which refined versions often exist. Our goal is rather to give a less constructive method than in former works, which results in much shorter proofs.

\section{Strategy and hypotheses for the convergence} 
Let us explain our idea to go rigorously  from \eqref{general} to \eqref{general_effective}. Given a smooth and compactly supported test field $\varphi$ on $\R^3$, the point is to  introduce a solution  $\phi_n = \phi_n[\varphi]$ of 
\begin{equation} \label{corrector}
 L^* \phi_n = -L^*_{e} \varphi \quad \text{in } \: \Omega_n 
\end{equation}
completed with appropriate boundary conditions at $\pa \Omega_n$. By appropriate, we mean that for natural extensions of $u_n$ and $\phi_n$ inside $\cup_i B_{i,n}$, still denoted $u_n$ and $\phi_n$,  one has the identity 
\begin{equation} \label{formal_identity}  
\int_{\Omega_n} (L u_n - g)   \, (\varphi -  \phi_n) = \int_{\R^3} (L u_n + L_e u_n - g) \varphi +  \int_{\R^3} g  \phi_n
\end{equation}
which in the case without boundaries would come from the formal calculations: 
$$  \int L u_n \, (\varphi -  \phi_n) = \int L u_n \varphi  -  \int u_n L^* \phi_n = \int L u_n \varphi  +  \int u_n L^*_{e} \varphi = \int (L + L_e) u_n \varphi. $$
By \eqref{general}, it will result from \eqref{formal_identity} that 
$$ \int_{\R^3} (L u_n + L_e u_n - g) \varphi +  \int_{\R^3} g  \phi_n = 0. $$
Hence, any accumulation point  $u$ of $u_n$ satisfies 
$$ |\int_{\R^3} (L u + L_e u - g) \varphi| \le  \limsup_n \Big| \int_{\R^3} g  \phi_n \Big|$$
To obtain \eqref{general_effective} in the limit $n \rightarrow +\infty$, or up to $o(\lambda)$ term,  it is then enough to show that $\phi_n$ converges to $0$ weakly, or up to a $o(\lambda)$ term.

\mspace
Let us remark that in the case of systems \eqref{Diff} and \eqref{Sto}, our solution $\phi_n$ is connected to the abstract test function $w^\eps$ that appears in  \cite{MR1493040} for the derivation of the strange term,  or in \cite{MR1079189}  for the derivation of the Brinkman term. We shall expand on this connection in Remark \ref{rem_1} below.  Nevertheless, in both \cite{MR1493040}  and \cite{MR1079189} , the concrete construction of  $w^\eps$ is quite tedious (and restricted to the periodic framework). Here, the direct introduction of $\phi_n$ simplifies the derivation:  its convergence to zero relies on relatively soft  arguments, in contrast with the use of correctors. 

\mspace
Namely, in the case where $L_e = L_e[\rho]$, the idea is that 
$$ \phi_n \approx \phi^1_{n}, \quad  \phi_n^1 :=  G^* \star \left( - L^*_e\bigl[\rho\bigr]\varphi\right)   -   G^* \star \left( - L^*_e\bigl[n^{-1} \sum \delta_{x_i}\bigr]\varphi\right) $$
with $G^*$  the fundamental solution of $L^*$. The first term in the definition of $\phi_n^1$ is to take care of the source term in \eqref{corrector}, while the second one is to take care of the boundary conditions at the spheres, neglecting the interaction between them. Note that $\phi_n^1$ goes formally to zero under the simple assumption \eqref{empirical}. Actually,  the main point, from which separation assumptions stem, is to show that  $\phi_n - \phi^1_n$ vanishes as $n \rightarrow +\infty$, or is $o(\lambda)$. Assumptions on the separation between the centers  stem from this requirement.  

\mspace
More precisely, in the context of system \eqref{Diff}, with $r_n = \frac{1}{n}$,  our assumptions are
\begin{equation} \label{H1} \tag{H1} 
 \exists c > 2, \quad  \forall  i \neq j, \quad  |x_i - x_j| \ge \frac{c}{n}, 
\end{equation}
\begin{equation} \label{H2}  \tag{H2}
\begin{aligned}
\limsup_{n \rightarrow +\infty}  \:   & \sum_i  \int_{B_i}  \hspace{-0.2cm}\bigg( n^2  \Big|  \frac{1}{n}\sum_{j \neq i} G(x-x_j) - \int_{\R^3} G(x-y) \rho(y) dy   \Big|^2    
+    \frac{1}{n^2}  \Big|  \sum_{j \neq i} \na G(x-x_j) \Big|^2  \bigg) dx = 0.
\end{aligned}
\end{equation}

\mspace
We have noted $x_i := x_{i,n}$, $B_i := B_{i,n}$ for short, and  $G(x) := \frac{1}{4\pi |x|}$ the kernel of $-\Delta$.  Settings in which \eqref{H1}-\eqref{H2} are satisfied, including periodic and random  i.i.d. cases, will be discussed in the last section \ref{secassumptions}. Let us already point out that \eqref{H1} and \eqref{H2} both express some kind of separation on the centers $x_{i,n}$. While \eqref{H1} is a mild requirement on their minimal distance, \eqref{H2} is a  condition of non-clustering in the mean, in the spirit of \eqref{weakseparation}. Let us note that if $G$ were smooth, \eqref{H2} would be a simple consequence of \eqref{empirical}. The point here is  the singularity of $G$ at the origin, which penalizes points $x_i, x_j$ close from one another. We state
\begin{theorem} \label{theo_Diff}
Let $g \in L^{\frac{6}{5}}(\R^3)$, $r_n = \frac{1}{n}$. Under \eqref{empirical}-\eqref{H1}-\eqref{H2}, the solution $u_n$  of \eqref{Diff} converges weakly in $\dot{H}^1(\R^3) \cap L^6(\R^3)$ to the solution $u$ of \eqref{Diff_eff}. 
\end{theorem}
\noindent
In the context of the Stokes equation \eqref{Sto}, assumption \eqref{H2} must be modified into 
\begin{equation} \label{H2'}  \tag{H2'}
\begin{aligned}
\limsup_{n \rightarrow +\infty}  \:    \sum_i  \int_{B_i} \bigg( n^2  \Big|  \frac{1}{n} & \sum_{j \neq i} G_{St}(x-x_j) - \int_{\R^3} G_{St}(x-y) \rho(y) dy   \Big|^2    
+    \frac{1}{n^2}  \Big|  \sum_{j \neq i} \na G_{St}(x-x_j) \Big|^2   \\
+    \frac{1}{n^4} \Big| & \sum_{j \neq i} R_{St}(x - x_j)  \Big|^2   + \frac{1}{n^6}  \Big|  \sum_{j \neq i} \na R_{St}(x - x_j)  \Big|^2 \bigg)dx = 0
\end{aligned}
\end{equation}
with  $G_{St}(x) = \frac{1}{8\pi} \left( \frac{I}{|x|} +  \frac{x \otimes x}{|x|^3}\right)$ the kernel of $-\mathbb{P} \Delta$, and  $R_{st}(x) = \frac{1}{8\pi} \left(  \frac{I}{3|x|^3} - \frac{x \otimes x}{|x|^5} \right)$. 
\begin{theorem} \label{theo_Stokes}
Let $g \in L^{\frac{6}{5}}(\R^3)^3$, $r_n = \frac{1}{n}$. Under \eqref{empirical}-\eqref{H1}-\eqref{H2'}, the solution $u_n$  of \eqref{Sto} converges weakly in $\dot{H}^1(\R^3)^3 \cap L^6(\R^3)^3$ to the solution $u$ of \eqref{Sto_eff}. 
\end{theorem}
\noindent
Eventually, system \eqref{Con} can be analyzed in the regime where  the volume fraction  $\lambda_n = \lambda$ is small but independent of $n$. We rely this time on  assumptions 
\begin{equation} \label{A1} \tag{A1} 
 \exists c > 2, \quad  \forall  i \neq j, \quad  |x_i - x_j| \ge c \, r_n, 
\end{equation}
and for all  smooth $\varphi$, 
\begin{equation} \label{A2}  \tag{A2}
 \limsup_{n \rightarrow +\infty} \:  \frac{1}{n^2}\sum_i \int_{B_i} \Big| \sum_{j \neq i} \na V(x-x_j) \na \varphi(x_j)\Bigr|^2 dx \le \eta(\lambda) \|\na\varphi\|_{L^\infty}, \quad \eta(\lambda) \xrightarrow[\lambda \rightarrow 0]{} 0, 
\end{equation}
where $V(x) = \na \frac{1}{4\pi |x|}$. This last assumption will be discussed further in Section \ref{secassumptions}. It will be notably shown to be implied by the usual condition \eqref{minimal_distance}. We shall prove
\begin{theorem} \label{theo_Con}
Let $g \in L^{\frac{6}{5}}(\R^3)$.  Assume $\lambda_n  = \lambda$ for all $n$.  Let $u_{n,\lambda}$ the solution of \eqref{Con}.  Under  \eqref{empirical}-\eqref{A1}-\eqref{A2}, any accumulation point $u_\lambda$ of $u_{n,\lambda}$ in  $\dot{H}^1(\R^3) \cap L^6(\R^3)$ solves 
\begin{equation} \label{sys:u}
- \div \big( (1+3\lambda|K_\rho| \rho)  \na u \big) = g + R_\lambda \quad \text{in $\R^3$}, 
 \end{equation}
with a remainder $R_\lambda$ satisfying 
 \begin{equation} \label{bound_R_lambda}
 \langle R_\lambda , \varphi \rangle  \le o(\lambda) \|\na\varphi\|_{L^\infty}, \quad \forall \varphi \in \dot{H}^1(\R^3) \cap L^6(\R^3) \cap \text{Lip}(\R^3).
 \end{equation}
\end{theorem}


\section{Proofs} \label{sec_proofs}
\subsection{Proof of Theorem \ref{theo_Diff}} \label{subsec_Diff}
We extend $u_{n}$ by zero inside the balls, and still denote $u_{n}$ this extension.  A simple energy estimate of \eqref{Diff} together with Sobolev imbedding yields 
$$ \|\na u_{n}\|^2_{L^2(\R^3)} \le \|g\|_{L^{\frac{6}{5}}(\R^3)} \|u_{n}\|_{L^6(\R^3)} \le C  \|g\|_{L^{\frac{6}{5}}(\R^3)} \|\na u_{n}\|_{L^2(\R^3)} $$
so that $u_{n}$ is bounded in $\dot{H}^1(\R^3) \cap L^6(\R^3)$. Up to  a subsequence in $n$, it has a weak limit $u$. 
Now,  given $\varphi \in C^\infty_c(\R^3)$, we introduce $\phi_n$ the solution of 
 \begin{equation}
 - \Delta \phi_n =  - 4 \pi \rho \varphi  \quad \text{ in } \Omega_n, \quad \phi_n\vert_{\cup_i B_i} = \varphi.
 \end{equation}  
Taking $\varphi - \phi_n$ as a test function in \eqref{Diff}, we get by Green's formula: 
$$ \int_{\R^3} \na u_{n} \cdot \na \varphi  +  4 \pi   \int_{\R^3} \rho u_{n} \varphi = \int_{\R^3} g \varphi  - \int_{\R^3} g \phi_n. $$
We have used here that $u_{n}$ and $\varphi - \phi_n$ vanish inside the balls to replace integrals over $\Omega_n$ by integrals over $\R^3$. Now,  if we prove that the last term at the right-hand side goes to zero, then $u \in \dot{H}^1(\R^3) \cap L^6(\R^3)$ will be a variational solution of \eqref{Diff_eff}. As this solution is unique, the full sequence $u_{n}$ will converge to $u$, proving the theorem. As $g$ is arbitrary in $L^{\frac{6}{5}}(\R^3)$, we need to prove that $\phi_n$ converges weakly  to zero in $L^6(\R^3)$. By a standard energy estimate on $\phi_n - \varphi$, which vanishes at the balls, it is easily seen that $\phi_n$ is bounded in $\dot{H}^1(\R^3) \cap L^6(\R^3)$, so that it is enough to show convergence of $\phi_n$ to zero in the sense of distributions.    

\mspace
Actually, we can further simplify the analysis by introducing the solution $\Phi_n$ of 
 \begin{equation} \label{equation_Phi_n}
 - \Delta \Phi_n =  - 4 \pi \rho   \quad \text{ in } \Omega_n, \quad \Phi_n\vert_{\cup_i B_i} = 1.   
 \end{equation}  
Indeed, let $R_n = \phi_n  - \Phi_n \varphi$,  that satisfies 
 \begin{equation}
 - \Delta R_n =   2 \na \Phi_n \cdot \na \varphi + \Phi_n \Delta \varphi  \quad \text{ in } \Omega_n, \quad R_n\vert_{\cup_i B_i} = 0.   
 \end{equation}  
If we prove that $\Phi_n \rightarrow 0$  strongly  in  $L^2_{loc}$, then, by a standard energy estimate, $R_n \rightarrow 0 $ strongly  in $H^1_{loc}$, and eventually  
$\phi_n \rightarrow 0$  strongly  in $L^2_{loc}$.

\mspace
For $\eta > 0$, we denote $\delta_n^\eta = \frac{1}{4\pi \eta^2 n} \sum_i s_\eta(\cdot - x_i)$ where $s_\eta$ is the surface measure on the sphere of radius $\eta$.  To prove strong  convergence of $\Phi_n$ to zero in $L^2_{loc}$, we split $\Phi_n = \Phi_n^1 + \Phi_n^2$, where 
$$ \Phi_n^1 =    4\pi \sum_i \mathcal{G}(n(x-x_i)) -  4\pi G \star \rho $$
with $G(x) = \frac{1}{4\pi|x|}$ the fundamental solution of $-\Delta$, $\mathcal{G}(x) = G(x)$ if  $|x| \ge 1$, $\mathcal{G}(x) = \frac{1}{4\pi}$  if  $|x| \le 1$. We compute
\begin{equation} \label{eq:Phi1n_Lap}
- \Delta \Phi_n^1 =  4 \pi  \delta_n^{1/n}   - 4 \pi \rho  \quad \text{ in } \R^3 
\end{equation}
and $\displaystyle \Phi_n^1(x) = 1 +   \frac{4\pi}{n} \sum_{j \neq i} G(x-x_j) -  4\pi \int_{\R^3} G(x-y) \rho(y) dy \: $ for $x \in B_i$, for all $i$. Hence,
\begin{equation}
- \Delta \Phi_n^2 = 0 \quad \text{ in } \Omega_n, \quad \Phi_n^2\vert_{B_i}(x) = -\frac{4\pi}{n} \sum_{j \neq i} G(x-x_j) +  4\pi \int_{\R^3} G(x-y) \rho(y) dy
\end{equation}

\mspace
Now, we use that under assumption \eqref{H1},  there is $C > 0$ such that for all $w = (w_i)_{1\le i \le n}$ in $\prod_i H^1(B_i)$, the solution $u[w]$ of 
$$ -\Delta u[w] = 0 \quad \text{ in } \: \Omega_n, \quad u[w]\vert_{B_i} = w_i  \text{ for all } i $$
satisfies 
\begin{equation} \label{estim_Lap}
\begin{aligned}
\|\na u[w]\|^2_{L^2(\R^3)} &  \le   C \sum_i \left( n^2 \|w_i\|_{L^2(B_i)}^2 + \|\na w_i\|^2_{L^2(B_i)} \right) 
\end{aligned}
\end{equation}
We refer to Lemma \ref{lemma_Sto} for a proof in the slightly more difficult Stokes case. Note that the factor $n^2$ at the right-hand side of the first inequality is consistent with scaling considerations. 
By applying this inequality with $w_i = -\frac{4\pi}{n}\sum_{j \neq i} G(x-x_j) +  4\pi \int_{\R^3} G(x-y)\rho(y) dy$, noticing that 
$\sum_i \int_{B_i} \big|  \int_{\R^3} \na G(x - y) \rho(y)dy   \big|^2  dx = O(\frac{1}{n^2})$,  and combining  with \eqref{H2},  we find that $\|\na \Phi^2_n\|_{L^2(\R^3)}$ goes to zero, so that  $\Phi^2_n$ goes to zero strongly in $H^1_{loc}$.   

\mspace
The last step is to show that $\Phi^1_n$ goes strongly to zero in $L^2_{loc}$. As the right-hand side of \eqref{eq:Phi1n_Lap} is bounded in $W^{-1,p}(\R^3)$ for any $p < \frac{3}{2}$, 
$\Phi^1_n$ is bounded in $W^{1,p}_{loc}$  for any $p < \frac{3}{2}$, with compact embedding in $L^q_{loc}$ for any $q < 3$. Hence, it is enough to prove convergence to zero in the sense of distributions. Let $\psi \in C^\infty_c(\R^3)$. Clearly the function $\Psi = \Delta^{-1} \psi$ is well-defined and  smooth. We get 
$$ \int_{\R^3}  \Phi^1_n \, \psi = \int_{\R^3} \Phi^1_n \Delta \Psi = \langle \Delta \Phi^1_n , \Psi \rangle  =  4 \pi  \langle \rho - \delta_n^{1/n} , \Psi \rangle \rightarrow 0$$
by the assumption \eqref{empirical}. This concludes the proof. 

\begin{remark} \label{rem_1}
Readers familiar with article \cite{MR1493040} may have recognized a connection  between the solution $\Phi_n$ of  \eqref{equation_Phi_n}, and the abstract couple $(w^\eps , \mu)$ involved in assumptions  (H.1)-(H.5) of \cite{MR1493040}.  Namely, setting $\eps = r_n$, one can check that $\Phi_n = 1 - w^\eps$ and $\mu  = 4\pi \rho$. However, the approach in the present paper remains quite distinct from the one in \cite{MR1493040}. In \cite{MR1493040}, the construction of the corrector $w^\eps$ is quite technical, and the identification of $\mu$ is done {\it a posteriori}, once $w^\eps$ has been built. Here, $\Phi_n$ (and the extra term  $4\pi \rho$) are introduced {\it a priori}, through equation \eqref{equation_Phi_n}. This allows for a much shorter derivation of the strange term, without any periodicity assumption.  
\end{remark}

\subsection{Proof of Theorem \ref{theo_Stokes}} 
Standard estimates show that the sequence $u_n$ of solutions of \eqref{Sto} is bounded in $\dot{H}^1(\R^3)^3 \cap L^6(\R^3)^3$. Let $\varphi \in C^\infty_c(\R^3)^3$ a divergence-free vector field. We introduce  the solution $\phi_n$ of 
\begin{equation}
-\Delta \phi_n + \na q_n = -6\pi \rho \varphi, \quad  \div \phi_n = 0 \text{ in } \Omega_n, \quad \phi_n\vert_{B_i} = \varphi. 
\end{equation}
Arguing exactly as in Paragraph \ref{subsec_Diff}, to show the weak convergence of $u_n$ to the solution $u$ of \eqref{Sto_eff}, it is enough to prove the convergence of $\phi_n$ to zero  in the sense of distributions. We consider this time the matrix-valued solution $\Phi_n$ of 
\begin{equation} \label{eq_cap_Phin}
-\Delta \Phi_n + \na Q_n = -6\pi \rho I, \quad  \div \Phi_n = 0 \text{ in } \Omega_n, \quad \Phi_n\vert_{B_i} = I,
\end{equation}
where $I$ denotes the identity matrix. With $R_n = \phi_n - \Phi_n \varphi$, and $S_n = q_n - Q_n \cdot \varphi$, we find
\begin{equation*}
\begin{aligned}
-\Delta R_n + \na S_n =  2 \pa_i \Phi_n \pa_i \varphi  + \Phi_n \Delta \varphi + (\na \varphi) Q_n, \quad \div R_n = \Phi_n : \na \varphi \quad \text{ in }\Omega_n, \quad R_n\vert_{B_i} = 0.  
\end{aligned}
\end{equation*}
We admit for the time being 
\begin{lemma} \label{lem:bogos}
Let $q > 3$, and $U$ a smooth bounded domain containing $K$. There exists a family of  operators  $\mB_n : L^q_0(\Omega_n \cap U) \rightarrow H^1_0(\Omega_n \cap U)$ 
such that $\div \mB_n h = h$  for all $h \in  L^q_0(\Omega_n \cap U)$ and such that $\|\mB_n\|_{\mathcal{L}(L^q_0, H^1_0)}$ is bounded uniformly in $n$. 
\end{lemma} 
\noindent
 Let $U$ containing both $K$ and the support of $\varphi$. One can check easily that  $\Phi_n : \na \varphi$  has zero average in both $U$ and $\Omega_n \cap U$. Extending 
 $\mB_n  (\Phi_n : \na \varphi)$ by zero outside $U$, the field $R'_n = R_n - \mB_n  (\Phi_n : \na \varphi)$  satisfies 
\begin{equation*}
\begin{aligned}
-\Delta R'_n + \na S'_n =  2 \pa_i \Phi_n \pa_i \varphi  + \Phi_n \Delta \varphi + (\na \varphi) Q_n + \Delta \mB_n  (\Phi_n : \na \varphi) \quad \text{in } \Omega_n, 
\end{aligned}
\end{equation*}
plus divergence-free and homogeneous Dirichlet conditions. By the previous lemma and a standard estimate, if we prove that  $\Phi_n \rightarrow 0$  strongly  in $L^q_{loc}$ for some $q > 3$ and that $Q_n \rightarrow 0$ strongly in $H^{-1}(U)$  (for some  appropriate extension in $\cup_i B_i$) , then $R'_n \rightarrow 0$ strongly  in  $H^1_{loc}$,  and so does $R_n$.  Eventually, $\phi_n$ will go to zero. 

\mspace
Therefore, we decompose $\Phi_n = \Phi^1_n + \Phi^2_n$, $Q_n = Q_n^1 + Q_n^2$. This time,  
\begin{align*}
\Phi^1_n(x) = 6\pi  \sum_{i} \mathcal{G}_{St}(n(x-x_i)) -  6\pi G_{St} \star \rho.   
\end{align*}
Here,  $G_{St}(x) =  \frac{1}{8\pi} \left( \frac{I}{|x|} +  \frac{x \otimes x}{|x|^3} \right)$ is the kernel of $-\mathbb{P} \Delta$, and  for   
$R_{St}(x) = \frac{1}{8\pi} \left(  \frac{I}{3|x|^3} -  \frac{x \otimes x}{|x|^5} \right)$, we have set
$$\mathcal{G}_{St}(x) = G_{St}(x) + R_{St}(x)   \quad  \text{if } \: |x| \ge 1, \quad  \mathcal{G}_{St}(x) = \frac{I}{6 \pi} \quad  \text{if } \:  |x| \le 1. $$ 
A tedious calculation shows that 
\begin{equation} \label{eq:Phi1n}
- \Delta \Phi_n^1 + \na Q^1_n =  6 \pi (\delta_n^{1/n}    -  \rho) I, \quad \div \Phi_n^1 \quad \text{ in } \R^3.  
\end{equation}
while 
\begin{align*}
& - \Delta \Phi^2_n  + \na Q_2^n = 0, \quad \div \Phi^2_n = 0 \quad \text{in } \: \Omega_n
\end{align*}
with boundary conditions 
\begin{align*}
 \Phi^2_n\vert_{B_i}(x) & =  -\frac{6\pi }{n}\sum_{j \neq i} G_{St}(x-x_j) + 6 \pi G_{St} \star \rho(x) -  \frac{6\pi }{n^3} \sum_{j \neq i} R_{St}(x-x_j) , \quad x \in B_i.
 \end{align*}
To estimate $\Phi^2_n$, we use the following (see below for a proof):
\begin{lemma} \label{lemma_Sto}
Under \eqref{H1}, inequality  \eqref{estim_Lap} is true for all (matrix or vector-valued) family $(w_i)_{1\le i \le n}$ in $\displaystyle \prod_i W^{1,\infty}(B_i)$ satisfying the compatibility condition $\int_{\pa B_i} w_i \cdot n = 0$ for all $i$, and  $u[w]$ the solution of 
$$ -\Delta u[w] + \na p[w] = 0, \quad \div u[w] = 0 \quad \text{ in } \: \Omega_n, \quad u[w]\vert_{B_i} = w_i  \text{ for all } i. $$
\end{lemma}
\noindent
Using this lemma and \eqref{H2'}, we find that  $\Phi^2_n$ converges strongly to zero in $\dot{H}^1(\R^3) \cap L^6(\R^3)$. Thanks to Lemma \ref{lem:bogos},  we can then show that for any $p < \frac{3}{2}$, the pressure term $Q^2_n$, normalized to be mean free in $\Omega_n \cap U$ and extended by zero in $\cup_i B_i$, converges strongly in $L^{p}(U)$. Indeed,  for any  $q > 3$ and any $h \in L^q(U)$, denoting $h_0 = h -  \dashint_{\Omega_n \cap U} h$, we find 
\begin{align*}
\int_{U} Q^2_n h & =   \int_{\Omega_n \cap U} Q^2_n h_0 =  \int_{\Omega_n \cap U} Q^2_n \div \mB_n h_0  = -\int_{\Omega_n \cap U} \na Q^2_n  \mB_n h_0 = - \int_{\Omega_n \cap U} \Delta \Phi^2_n  \mB_n h_0 \\
& =   \int_{\Omega_n \cap U} \na \psi^2_n \cdot  \na  \mB_n h_0 \le C \|\na  \Phi^2_n\|_{L^2(\Omega_n \cap U)}  \|\na \mB_n h_0\|_{L^2(\Omega_n \cap U )} \\
 & \le  C \|\na  \Phi^2_n\|_{L^2(\R^3)}  \|h_0\|_{L^q(U)} \le   C' \|\na  \Phi^2_n\|_{L^2(\R^3)}  \|h\|_{L^q(U)}
\end{align*}
The strong convergence of $Q^2_n$ in $L^{p}$, $p=q'$,  follows. 

\mspace
We finally have to look at the convergence of $(\Phi^1_n, Q^1_n)$ solving \eqref{eq:Phi1n}. The source  $6 \pi (\delta_n^{1/n}    -  \rho) I$ is compactly supported, and converges to zero in the sense of measures, in particular weakly in $W^{-1,p}(\R^3)$ for any $p < \frac{3}{2}$.  Moreover, for an appropriate normalization of the pressure $Q^1_n$, we have the estimate 
$$ \|\na \Phi^1_n\|_{L^p(\R^3)} + \|\Phi^1_n\|_{L^{3p/(3-p)}(\R^3)} + \|Q^1_n\|_{L^p(\R^3)} \le C \|6 \pi (\delta_n^{1/n}    -  \rho) I\|_{W^{-1,p}(\R^3)} $$
see \cite[Theorem 4.2.2]{MR1284205}.  This implies that for any $p < \frac{3}{2}$, $\Phi^1_n \rightarrow 0$ weakly in  $L^{\frac{3p}{3-p}}(\R^3)$ and strongly in $L^q_{loc}$ for any $q < \frac{3p}{3-p}$, while  $Q^1_n \rightarrow 0$ weakly in $L^p(\R^3)$. 

\mspace
Collecting all previous facts, we see that $\Phi_n \rightarrow 0$  strongly in $L^q_{loc}$ for any $q < 3$, and $Q_n \rightarrow 0$ weakly in $L^p(U)$ for any $p < \frac{3}{2}$, therefore strongly in $H^{-1}(U)$. To gain on the exponent $q$ for $\Phi_n$ and conclude the proof, it is enough to show that $\na \Phi_n$ is bounded in $L^2(\R^3)$, which implies strong convergence in $L^q_{loc}$ for any $q < 6$.  Such uniform bound can be obtained easily: by use of standard Bogovski operator, one can find a compactly supported (matrix-valued) $H^1$ field $\Phi_I$ which is divergence-free  and such that $\Phi_I = I$ in a neighborhood of $K \supset \cup_i B_i$. A uniform $L^2$ bound on $\na (\Phi_n - \Phi_I)$ then follows from a standard energy estimate, as $\Phi_n - \Phi_I$ is zero in the balls.

\mspace
{\bf Proof of Lemma \ref{lem:bogos}}. Let $q > 3$, and $h \in L^q_0(\Omega_n \cap U)$. We extend $h$ by zero in $\cup_i B_i$.  First, we introduce a standard Bogovski operator $\mB : L^q(U) \rightarrow W^{1,q}_0(U)$ which is continuous and satisfies $\div \mB h = h$. The next step is to find a  field $w_n = w_n[h]$ in $H^1_0(U)$  such that 
$$ \div w_n = 0 \: \text{ in } U \cap \Omega_n, \quad w_n\vert_{B_i} =   \mB h\vert_{B_i}, \quad \|w_n\|_{H^1_0} \le C \|\mB h\|_{W^{1,q}_0}, \quad C \text{ independent of $n$}.   $$
Indeed, $\mB_n h = \mB h - w_n$ will have the required properties. As the balls $nB_i =  B(n x_i, 1)$ have radius $1$, application of  usual Bogovski operators  yield for each $1 \le i \le n$ a divergence-free function $W_i \in  H^1_0(B(nx_i, 2))$  such that 
$$W_i\vert_{n B_i} = \mB h(n^{-1} \cdot)\vert_{n B_i}, \quad \|W_i\|_{H^1_0(B(nx_i,2))} \le C \| \mB h(n^{-1}\cdot)\|_{H^1(n B_i)}, \quad C \: \text{independent of $n$}. $$
We set $w_n(x) = \sum_i W_i(nx)$. Clearly, $w_n$ is divergence-free, and  belongs to $H^1_0(U)$. Moreover,  assumption \eqref{H1} implies that the balls $B(nx_i, 2)$ are disjoint, so that  $w_n\vert_{B_i} = W_i(n \cdot)\vert_{B_i} = \mB h\vert_{B_i}$, and 
\begin{align*}
 \|w_n\|^2_{H^1_0(U)}  & \le C  \|\na w_n\|^2_{L^2(U)} =  \sum_{i} \|\na (W_i(n \cdot))\|_{L^2(B(x_i, 2/n))}^2 \\ 
 & \le C \sum_i \big( n^2 \|\mB h\|_{L^2(B_i)}^2 +  \|\na \mB h\|_{L^2(B_i)}^2 \bigr) \\
 & \le \frac{4\pi C}{3} \|\mB h\|_{L^\infty(U)} + C   \|\na \mB h\|_{L^2(U)}^2 \le C' \|\mB h\|_{W_0^{1,q}(U)}.  
  \end{align*}
where the last inequality involves the Sobolev embedding $W^{1,q} \hookrightarrow L^\infty$. This ends the proof. 

\mspace
{\bf Proof of Lemma \ref{lemma_Sto}}. By the classical variational characterization  of the Stokes solution,  
$$ \|\na u[w]\|_{L^2(\R^3)}^2 \le \|\na v[w]\|_{L^2(\R^3)}^2 $$
for any divergence-free vector field $v = v[w]$ such that $v\vert_{B_i} = w_i$ for all $i$. Hence, it is enough to prove that there exists such a $v$ satisfying the bound \eqref{estim_Lap}. One proceeds as in the proof of Lemma \ref{lem:bogos}: one looks for $v$ under the form $v = \sum_i W_i(nx)$, where $W_i \in  H^1_0(B(nx_i, 2))$ is provided through the use of Bogovski operator:  
$$W_i\vert_{n B_i} = w_i(n^{-1} \cdot)\vert_{n B_i}, \quad \|W_i\|_{H^1_0(B(nx_i,2))} \le C \| w_i(n^{-1}\cdot)\|_{H^1(n B_i)}, \quad C \: \text{independent of $n$}. $$
Estimate \eqref{estim_Lap} follows as in the proof of Lemma \ref{lem:bogos}. 

\begin{remark}
One can extend the result for the Stokes system \eqref{Sto} to the Navier-Stokes system 
\begin{equation} \label{NSto}
 - \Delta u_n  + u_n \cdot \na u_n + \na p_n = g, \quad \div u_n = 0  \quad  \text{in} \:  \Omega_n, \quad
 u_n\vert_{B_{i,n}} = 0, \quad 1 \le i \le n. 
 \end{equation}
In such a case, we can not use directly the test function $\varphi - \phi_n$, as $\Omega_n$ is unbounded : for instance, the  extra  term $ \int_{\R^3} (u_n \cdot \na u_n) \cdot  \phi_n = - \int_{\R^3} (u_n \otimes u_n) : \na \phi_n$ due to the nonlinearity is not {\it a priori} properly defined. To circumvent this problem, one can use as a test function  $\chi (\varphi - \phi_n) = \varphi - \chi \phi_n$, where $\chi \in C^\infty_c(\R^3)$ satisfies $\chi = 1$ in a large open ball containing  all the balls $B_i$ and the support of $\varphi$. Compared to the Stokes case, there are then several extra integrals in the variational formulation: 
\begin{align*}
 &I_n := \int_{\R^3} (u_n \otimes  u_n) : \na (\varphi - \chi \phi_n), \quad  J_n := \int_{\R^3}  \na u_n : (\na \chi \otimes \phi_n) - (u_n \otimes \na \chi) : \na \phi_n \\
 &K_n := \int_{\R^3} p_n \na \chi \cdot \phi_n,  \quad L_n :=  \int_{\R^3} q_n \na \chi \cdot u_n, \quad  M_n :=  \int_{\R^3} (1-\chi) g \phi_n
 \end{align*}
By a simple energy estimate and Sobolev embedding, still valid in the Navier-Stokes case,  $u_n$ converges to $u$   weakly in  $\dot{H}^1 \cap L^6$, thus strongly in $L^p_{loc}$ for any $p < 6$. From the proof above, we also have that $\phi_n$ converges to zero  weakly in  $\dot{H}^1 \cap L^6$, thus strongly in $L^p_{loc}$ for any $p < 6$. From there, clearly: $ I_n \rightarrow \int_{\R^3} (u \otimes  u) : \na \varphi, \quad J_n \rightarrow 0, \quad M_n \rightarrow 0$.  As regards $K_n$, we  use the fact that we can choose an open  neighborhood $U$ of the support of $\na \chi$,  independent of $n$, such that 
$$ - \Delta u_n  + \na p_n = - \div(u_n \otimes u_n) + g \quad \text{ in } U. $$
Choosing $p_n$ so that it has zero average in $U$, we deduce from standard estimates that
$$ \|p_n\|_{L^2(U)} \le C_U (\|\na u_n\|_{L^2(U)} + \|u_n \otimes u_n \|_{L^2(U)} + \|g\|_{L^{6/5}(U)}  \le C'. $$   
Using the strong convergence of $\phi_n$ in $L^2_{loc}$  and this uniform bound, we get $K_n \rightarrow 0$. By inverting the roles of $(u_n,p_n)$ and $(\phi_n,q_n)$, we get $L_n \rightarrow 0$ as well, resulting in the limit system
\begin{equation} \label{NSto_eff}  
-\Delta u+ u \cdot \na u + \na p + 6 \pi \rho u = g, \quad \div u = 0 \text{ in } \R^3. 
\end{equation}
\end{remark}

\begin{remark}
The same methodology applies to the system considered in \cite{MR2398959}-\cite{MR3851058}, where the homogeneous Dirichlet condition in \eqref{Sto} is replaced by an inhomogeneous one:
$$\forall 1 \le i \le n, \quad u_n\vert_{B_{i,n}} = V_{i,n}$$
 for a collection of constant velocities $(V_{i,n})_{1 \le i \le n}$ satisfying 
$$ \sup_n  \frac{1}{n} \sum_i |V_{i,n}|^2 < +\infty, \quad j_n \, :=  \, \frac{1}{n} \sum V_{i,n} \delta_{x_{i,n}} \xrightarrow[n \rightarrow \infty]{} j \quad \text{ in the sense of measures}.  $$
To handle this situation, there is one extra step: to show that the solution $v_n$ of 
\begin{equation*} 
 -\Delta v_n + \na q_n = -6 \pi j, \quad \div v_n = 0 \quad \text{ in } \Omega_n, \quad v_n\vert_{B_{i,n}} = V_{i,n} 
 \end{equation*}
 converges to zero weakly in $\dot{H}^1(\R^3)^3 \cap L^6(\R^3)^3$. Assuming such convergence, one can then notice that $u_n - v_n$ satisfies a system of type \eqref{Sto}, with source term $g + 6\pi j$ replacing $j$. By Theorem \ref{theo_Stokes}, the weak limit of $u_n$, which coincides with the weak limit of $u_n - v_n$, satisfies 
 $$ -\Delta u + \na p + 6\pi u = g + 6\pi j, \quad \div u \quad \text{ in } \: \R^3. $$
To prove the convergence of $v_n$ to zero,  we mimic what was done for $\Phi_n$, solution of \eqref{eq_cap_Phin}. Just as we decomposed $\Phi_n = \Phi^1_n + \Phi^2_n$, we write $v_n = v^1_n + v^2_n$, with approximate solution  
$$ v^1_n  := 6\pi  \sum_{i} \mathcal{G}_{St}(n(x-x_i)) V_{i,n} -  6\pi G_{St} \star j.  $$
 Weak convergence of $v^1_n$ to zero comes from  the convergence of $j_n$ to $j$, while the remainder $v_2^n$ goes strongly to zero in $H^1_{loc}$ under the assumption 
\begin{equation*}
 \begin{aligned}
\limsup_{n \rightarrow +\infty}  \:    \sum_i  \int_{B_i} \bigg( n^2  \Big|  \frac{1}{n} & \sum_{j \neq i} G_{St}(x-x_j)V_{j,n} - \int_{\R^3} G_{St}(x-y) j(y) dy   \Big|^2    
+    \frac{1}{n^2}  \Big|  \sum_{j \neq i} \na G_{St}(x-x_j) V_{j,n} \Big|^2   \\
+    \frac{1}{n^4} \Big| & \sum_{j \neq i} R_{St}(x - x_j) V_{j,n} \Big|^2   + \frac{1}{n^6}  \Big|  \sum_{j \neq i} \na R_{St}(x - x_j) V_{j,n} \Big|^2 \bigg)dx = 0
\end{aligned}
\end{equation*}
which is  similar to \eqref{H2'}, and holds in classical settings, see Section \ref{secassumptions} for discussion of \eqref{H2'}. 
\end{remark}

\subsection{Proof of Theorem \ref{theo_Con}}
We remind that  here, $r_n$ is of order $n^{-1/3}$, namely with a volume fraction $\lambda = n \frac{4}{3} \pi r_n^3$ positive and independent of $n$. 
As in other examples, the energy estimate and the Sobolev imbedding give that the sequence $(u_{n,\lambda})_{n \in \N}$ is bounded in  $\dot{H}^1(\R^3) \cap L^6(\R^3)$. After extraction, it converges weakly to some $u_\lambda$. Let $\varphi \in \dot{H}^1(\R^3) \cap L^6(\R^3) \cap C^\infty(\R^3)$. We consider $\phi_{n,\lambda}$ satisfying 
\begin{equation} \label{eq_phi_Con}
- \Delta \phi_{n,\lambda} = \div(3 \lambda |K_\rho|  \rho \na \varphi)  \: \text{ in } \Omega_n,   \: \int_{\pa B_i} \pa_\nu \phi_{n,\lambda} = - \int_{\pa B_i} 3 \lambda |K_\rho| \rho \pa_\nu \varphi,  \: \phi_{n,\lambda}\vert_{B_i}= \varphi\vert_{B_i} + \text{ cst},  \forall i. 
\end{equation}
Testing $\varphi - \phi_{n,\lambda}$ in  \eqref{Con}, we obtain after a few integrations by parts 
$$ \int_{\R^3} \na u_{n,\lambda} \cdot (1 + 3 \lambda |K_\rho| \rho)  \na \varphi  = \int_{\R^3} g \varphi - \int_{\R^3} g \phi_{n,\lambda}.  $$
As $n$ goes to infinity, the left-hand side converges to $\int_{\R^3} \na u_\lambda \cdot (1 + 3 \lambda f)  \na \varphi$, which implies that the last term at the right-hand side, that can be written in the abstract form $\langle R_{n, \lambda} , \varphi \rangle$ for an element $R_{n, \lambda}$ in the dual of $\dot{H}^1(\R^3) \cap L^6(\R^3)$, converges to some  $\langle R_{\lambda} , \varphi \rangle$, with $R_\lambda$ in this same dual. 
To prove the theorem, it is enough to show that 
\begin{equation} \label{tobeproved}
\limsup_{n \rightarrow +\infty} |\int_{\R^3} g \phi_{n,\lambda}|  \le \lambda \eta(\lambda) \|\na\varphi\|_{L^\infty}, \quad  \eta(\lambda) \xrightarrow[\lambda \rightarrow 0]{} 0.
\end{equation}
Let us stress that restricting  to smooth test functions $\varphi$ instead of Lipschitz is no problem: indeed,  if \eqref{tobeproved} holds for smooth functions, we can apply it to $\rho_\delta \star \varphi$ with $\rho_\delta$ a mollifier. We deduce
$$ \langle R_\lambda, \rho_\delta \star  \varphi \rangle \le  \lambda \eta(\lambda) \|\na \rho_\delta \star \varphi \|_{L^\infty} \le  \lambda \eta(\lambda) \|\na \varphi \|_{L^\infty}  $$
As $\rho_\delta \star  \varphi$ converges in $\dot{H}^1(\R^3) \cap L^6(\R^3)$, the l.h.s. converges to $\langle R_\lambda , \varphi \rangle$, which yields \eqref{bound_R_lambda}.  

\mspace 
We split $\phi_{n,\lambda} = \phi^1_n + \phi^2_n$, with  
$$ \phi^1_n =        3 \lambda |K_\rho|  V \star (\rho  \na \varphi) - 4\pi r_n \sum_i \na \varphi(x_i) \cdot \mathcal{V}\left(\frac{x-x_i}{r_n}\right), $$
where $V(x) = \na \frac{1}{4\pi |x|}$, $\mathcal{V}(x) = V(x)$ for $|x| \ge 1$, $\mathcal{V}(x) = -\frac{x}{4\pi}$ for $|x| \le 1$. 

\mspace
We compute $ \Delta \mathcal{V} = \frac{3}{4\pi} x s^1$, with $s^\eta$ the surface measure on the sphere of radius $\eta$. Hence,  
\begin{equation} \label{eq_Phi1n}
\begin{aligned}
-\Delta \phi^1_n & = \div \Big(3 \lambda |K_\rho| \rho \na \varphi \Big)  +  \frac{3}{r_n}\sum_i \na \varphi(x_i) \cdot (x-x_i) s^{r_n}(x-x_i)     \\
& =  \div \Big(3 \lambda |K_\rho| \rho \na \varphi \Big)   -  3  \sum_i \na \varphi(x_i) \cdot \na 1_{B_i}   \quad \text{ in } \R^3  \\
&= \div \Big( 3 \lambda |K_\rho| \rho \na \varphi  - 3 \lambda |K_\rho|\Big( \frac{1}{n} \sum_i \frac{1}{|B_i|}1_{B_i} \na \varphi(x_i) \Big) \Big) \quad \text{ in } \R^3.
\end{aligned}
\end{equation} 
One also checks that $\int_{\pa B_i} \pa_\nu \phi^1_n = - \int_{\pa B_i} 3 \lambda \rho \pa_\nu \varphi$. We further decompose $\phi^2_n = \psi^2_n + \tilde \psi^2_n$, with 
$$ - \Delta \psi^2_n = - \Delta \tilde \psi^2_n = 0, \quad \int_{\pa B_i} \pa_\nu \psi^2_n =   \int_{\pa B_i} \pa_\nu \tilde \psi^2_n = 0, $$
 and for all $i$: 
 \begin{align*}
 \psi^2_n \vert_{B_i}(x) & =   -3 \lambda |K_\rho|  V \star (\rho  \na \varphi)(x)  + 4\pi r_n \sum_{j \neq i} \na \varphi(x_j) \cdot V\left(\frac{x-x_j}{r_n}\right)  + cst \\
\tilde  \psi^2_n \vert_{B_i}(x) & = \varphi(x) - \na \varphi(x_i) \cdot (x-x_i) + cst.
 \end{align*}
By classical variational characterization, the function $\psi^2_n$, resp. $\tilde  \psi^2_n$, minimizes the Dirichlet integral among all functions $\psi \in \dot{H}^1(\R^3)$ satisfying 
$\psi\vert_{B_i} = \psi^2_n\vert_{B_i} + c_i$, resp. $\psi\vert_{B_i} = \tilde \psi^2_n\vert_{B_i} + c_i$,  on $B_i$, for some $c_i \in \R$, for all $i$. From there, one can proceed as in the proof of \eqref{estim_Lap}, {\it cf.} Lemma \ref{lemma_Sto}, and show that under \eqref{A1}: 
\begin{equation} \label{trace_control}
 \|\na \psi^2_n\|^2_{L^2(\R^3)} \le C  \sum_{i} \int_{B_i} |\na \psi^2_n|^2, \quad   \|\na \tilde \psi^2_n\|^2_{L^2(\R^3)} \le C  \sum_{i}  \int_{B_i} |\na \tilde \psi^2_n|^2. 
 \end{equation}
 note that, in contrast to \eqref{estim_Lap}, the right-hand side involves only gradients, consistently with the additional degree of freedom provided by possible addition of constants $c_i$ on $B_i$.  
In the case of $\tilde \psi^2_n$, we get 
\begin{equation}  \label{estimate_tildePsi2}
 \|\na \tilde \psi^2_n\|^2_{L^2(\R^3)} \le C \|\na^2 \varphi\|_{L^\infty} \sum_i \int_{B_i} |x-x_i|^2 dx = O(n^{-2/3}).
 \end{equation}
As regards $\psi^2_n$, we compute 
\begin{align*}
 \na \psi^2_n\vert_{B_i}(x) 
 &  =   -3\lambda |K_\rho|  \int   \rho(y) \na V(x-y) \na \varphi(y) dy +\frac{3\lambda |K_\rho|}{n}   \sum_{j \neq i} \na V(x-x_j)\na \varphi(x_j) \\
 &  = -b(x) + a_i(x), 
\end{align*}
where $\na V = -\frac{1}{4\pi}\na \frac{x}{|x|^3}$ is a Calderon-Zygmund kernel. Hence,   $b$ satisfies for all $q > 1$: 
$$\|b\|_{L^q(\R^3)} \le C_q \lambda \| \rho\na \varphi\|_{L^q(\R^3)} \le C_q \lambda  \| \rho\na \varphi\|_{L^q(K)}  $$
It follows by H\"older inequality that for all $q > 2$
\begin{equation} \label{ineqb}
 \sum_i \|b\|_{L^2(B_i)}^2 = \|b\|_{L^2(\cup B_i)}^2  \le \lambda^\frac{q-2}{q} \|b\|_{L^q(\R^3)}^{2} \le C \lambda^{3-\frac{2}{q}} \| \na \varphi\|_{L^{\infty}}^2. 
\end{equation}
Eventually, by \eqref{trace_control}, we find for all $\eps > 0$:
\begin{equation*} 
\begin{aligned}
\limsup_{n \rightarrow +\infty}  \|\na \psi^2_n\|^2_{L^2(\R^3)} & \le C_\eps \lambda^{3-\eps} \| \na \varphi\|_{L^{\infty}}^2  + \limsup_{n \rightarrow +\infty} \sum_i \|a_i\|_{L^2(B_i)}^2 \\
& \le (C_\eps \lambda^{3-\eps} + \eta(\lambda) \lambda^2) \|\na \varphi\|_{L^\infty(\R^3)}^2, \quad \eta(\lambda) \xrightarrow[\lambda \rightarrow 0]{} 0 
\end{aligned}
\end{equation*}
where the last inequality comes from assumption \eqref{A2}.  Together with \eqref{estimate_tildePsi2}, we deduce 
$$ \limsup_{n \rightarrow +\infty} \|\na \phi^2_n\|_{L^2(\R^3)} \le \eta(\lambda) \lambda \|\na \varphi\|_{L^\infty}, \quad \eta(\lambda) \xrightarrow[\lambda \rightarrow 0]{} 0. $$
To show \eqref{tobeproved} and conclude in this way the proof of the theorem, it is enough to show that $\phi^1_n$ goes weakly to zero as $n \rightarrow +\infty$ in $L^6(\R^3)$. Standard estimates on $\phi - \phi_{n, \lambda}$ show that $\phi_{n, \lambda}$, and from there $\phi_n^1$, is bounded uniformly in $n$ in $\dot{H}^1(\R^3) \cap L^6(\R^3)$, so that convergence to zero  in the sense of distributions is enough. It follows from a duality argument: for $h$ any smooth and compactly supported function, we introduce $H$ the solution of $\Delta H = h$ in $\R^3$, and thanks to \eqref{eq_Phi1n}, we write
\begin{align*}
\langle \phi^1_n , h \rangle & =  \langle \Delta \phi^1_n , H \rangle  =  3 \lambda |K_\rho| \,  \langle \rho\na \varphi  -   \frac{1}{n} \sum \frac{1}{|B_i|} 1_{B_i}\na \varphi(x_i), \na H \rangle 
\end{align*}
which is easily seen to go to zero as $n \rightarrow +\infty$, by \eqref{empirical}. 

\section{Further discussion of the hypotheses} \label{secassumptions}
We conclude this paper with further comments on the assumptions of our theorems. First, we show that the assumptions \eqref{H1} and \eqref{H2} in Theorem \ref{theo_Diff} are implied by the conditions
\begin{equation} \label{H1sharp} \tag{H1$^\sharp$}
\lim_{n \rightarrow +\infty} n  d_n = +\infty, \quad d_n = \inf_{i \neq j} |x_i - x_j|
\end{equation}
\begin{equation} \label{H2sharp} \tag{H2$^\sharp$}
\lim_{n \rightarrow +\infty} \frac{1}{n}   \sum_{i}  \Big( \frac{1}{n}\sum_{j \neq i} \frac{1}{|x_i - x_j|} - \int_{\R^3} \frac{\rho(y)}{|x_i -y|} dy\Big)^2 \: = 0 
\end{equation}
Clearly,  \eqref{H1sharp} is  stronger than \eqref{H1}. Then,  \eqref{H2} follows from  \eqref{H2sharp} if we prove that 
\begin{equation*}
\begin{aligned}
 \lim_n  & \:  \sum_i   \int_{B_i}\bigg(  \Big| \sum_{j \neq i} (G(x-x_j) -  G(x_i-x_j))     \Big|^2  \\
 & \:  +  n^2 \Big| \int_{\R^3} G(x-y) \rho(y) dy - \int_{\R^3} G(x_i-y) \rho(y) dy \Big|^2   +  \Big|  \frac{1}{n}\sum_{j \neq i} \na G(x-x_j) \Big|^2  \bigg) dx  = 0. 
\end{aligned}
\end{equation*}
We focus on the first term, as the second one is simpler and the third one  similar. We write: 
\begin{align*}
 & \sum_i   \int_{B_i} \Big| \sum_{j \neq i} (G(x-x_j) -  G(x_i-x_j))     \Big|^2
\\  \le \:  C  & \sum_i \int_{B_i}   \Big| \sum_{j \neq i}  \sup_{z \in [x_i, x]} |\na G(z-x_j)| |x-x_i|     \Big|^2 \\
\le \:   \frac{C'}{n^5}  &  \sum_i \big| \sum_{j \neq i} \frac{1}{|x_i -x_j|^2} \Big|^2 \: \le \:   \frac{C'}{(n d_n)^2}  \frac{1}{n} \sum_i \Big| \frac{1}{n}  \sum_{j \neq i} \frac{1}{|x_i -x_j|} \Big|^2
\: \le \:  \frac{C''}{(n d_n)^2} 
\end{align*}
where we deduced from  \eqref{H2sharp} that $\frac{1}{n} \sum_i \Big| \frac{1}{n}  \sum_{j \neq i} \frac{1}{|x_i -x_j|} \Big|^2$ is bounded uniformly in $n$. Convergence to zero follows then from \eqref{H1sharp}. 

\mspace
Note that \eqref{H2sharp} is in the same spirit as the  hypothesis \eqref{weakseparation} used in \cite{MR609184}, although a bit stronger. Both \eqref{H1sharp} and  \eqref{H2sharp} apply to many contexts. First, they are satisfied when the points are well-separated, {\it cf.} \eqref{minimal_distance}. Let us deduce \eqref{H2sharp} from \eqref{minimal_distance}. For all $j$, we set $\mathcal{B}_j = B\big(y_j, \frac{c}{2} n^{-1/3}\big)$, where $c$ is the constant appearing in \eqref{minimal_distance}. As $\frac{1}{|x|}$ is harmonic, we apply the mean value formula to write for any $i$: 
\begin{align*}
\frac{1}{n}\sum_{j \neq i} \frac{1}{|x_i - x_j|} \: & = \:    \frac{1}{n}\sum_{j \neq i} \frac{1}{|\mathcal{B}_j|} \int_{\mathcal{B}_j}  \frac{1}{|x_i - y|} dy  \: = \:  \frac{1}{n}\sum_{j}  \frac{1}{|\mathcal{B}_j|} \int_{\mathcal{B}_j}  \frac{1}{|x_i - y|} dy +  O(n^{-2/3})  \\
& =  \: \int_{\R^3} \frac{1}{|x_i- y|} g_n(y) dy   + O(n^{-2/3})
\end{align*}
where $\rho_n(y) =  \frac{1}{n}\sum_{j}  \frac{1}{|\mathcal{B}_j|}  1_{B_j}(y)$. By  \eqref{H2sharp}, the sequence $(\rho_n)_{n \in \N}$ is  bounded  in $L^\infty$, and by \eqref{empirical} is easily seen to converge weakly * to $\rho$. Hence, for all $x \in K$,  $f_n(x) =   \int_{\R^3} \frac{1}{|x- y|} \rho_n(y) dy$ converges to  $f(x) =   \int_{\R^3} \frac{1}{|x- y|} \rho(y) dy$. It is also easily seen that $(f_n)$ is equicontinuous over $K$, from which we deduce that 
$\sup_K   |f_n - f| \rightarrow 0$, and from there  \eqref{H2sharp}.

\mspace
The convergences in  \eqref{H1sharp} and \eqref{H2sharp} hold also in probability when the $x_i$ are random i.i.d variables with law $\rho$. More precisely,  it is well-known that in such setting, $n^\alpha d_n \rightarrow +\infty$ in probability for any $\alpha > \frac{2}{3}$, which is of course stronger than  \eqref{H1sharp}.  Moreover, 
$$ \lim_{n \rightarrow +\infty} \e \frac{1}{n}   \sum_{i}  \Big( \frac{1}{n}\sum_{j \neq i} \frac{1}{|x_i - x_j|} - \int_{\R^3} \frac{\rho(y)}{|x_i -y|} dy\Big)^2 \: = 0 $$ 
which implies again convergence in probability. To establish this last limit is very classical: we write 
\begin{align*}
   & \e \Big( \frac{1}{n}\sum_{j \neq i} \frac{1}{|x_i - x_j|} - \int_{\R^3} \frac{\rho(y)}{|x_i -y|} dy\Big)^2 = \e  \Big( \frac{1}{n}\sum_{j \neq i} \Big( \frac{1}{|x_i - x_j|} - \int_{\R^3} \frac{\rho(y)}{|x_i -y|} dy \Bigr) \Big)^2 + O(n^{-2}) \\
   & =  \frac{1}{n^2} \sum_{j,j' \neq i}  \e X_{ij} X_{ij'} + O(n^{-2}) = \frac{1}{n^2} \sum_{j \neq j' \neq i}  \e X_{ij} X_{ij'} + O(n^{-1})  
 \end{align*}
where $X_{ij} :=   \frac{1}{|x_i - x_j|} - \int_{\R^3} \frac{\rho(y)}{|x_i - y|} dx$ has mean zero and variance  
$$\e X_{ij}^2 = \int  \rho(x) \rho(z)  \Big|  \frac{1}{|x-z|} - \int \frac{\rho(y)}{|x-y|} dy \Big|^2 dx dz.$$ 
Eventually, for $j \neq j' \neq i$, we compute using the independance of the $x_k$'s: 
\begin{align*}
  \e X_{ij} X_{ij'}  & = \int_{(\R^3)^3} \rho(z) \rho(x) \rho(x')  \Big(\frac{1}{|z -x|} - \int_{\R^3} \frac{\rho(y)}{|z -y|} dy \Big) \, \Big(\frac{1}{|z -x'|}  - \int \frac{\rho(y')}{|z -y'|} dy' \Big) dz dx dx' \\
& =   \int_{\R^3} \rho(z)  \Big( \int_{\R^3} \frac{\rho(x)}{|z -x|} dx - \int_{\R^3} \frac{\rho(y)}{|z -y|} dy \Big)  \Big( \int_{\R^3} \frac{\rho(x')}{|z -x'|} dx' - \int_{\R^3} \frac{\rho(y')}{|z -y'|} dy' \Big) dz  = 0.  
\end{align*}
Hence, \eqref{H1sharp}-\eqref{H2sharp} holds in probability when the $x_i$ are i.i.d. random variables. We remind that a sequence is converging in probability if and only if any subsequence  has itself a subsequence that converges almost surely. Using this characterization, and applying the proof of Theorem \ref{theo_Diff}, we find that $u_n$ converges to $u$ solution of \eqref{Diff_eff} in probability, for any distance metrizing the weak  topology  of the ball of $\dot{H}^1 \cap L^6$ to which  all $u_n$ belong. 

\mspace
Similar considerations apply to the Stokes case. We leave to the reader to check that assumptions \eqref{H1}-\eqref{H2'} are implied by \eqref{H1sharp}- \eqref{H2'sharp}, with 
\begin{equation} \label{H2'sharp} \tag{H2'$^\sharp$}
\lim_{n \rightarrow +\infty} \frac{1}{n}   \sum_{i}  \Big( \frac{1}{n}\sum_{j \neq i} \frac{(x_i - x_j) \otimes (x_i - x_j)}{|x_i - x_j|^3} - \int_{\R^3} \frac{(x_i - y) \otimes (x_i - y)}{|x_i -y|^3} \rho(y) dy\Big)^2 \: = 0.  
\end{equation}
These stronger hypotheses are again verified under the assumption \eqref{minimal_distance}, or in the i.i.d. case. To this respect, the only real change is in showing that  \eqref{minimal_distance} implies \eqref{H2'sharp}. We use this time that $\frac{x}{|x|^3} = -\na \frac{1}{|x|}$ is harmonic and obtain 
\begin{align*} 
\frac{1}{n}\sum_{j \neq i} \frac{(x_i - x_j) \otimes (x_i - x_j)}{|x_i - x_j|^3} & =  \frac{1}{n}\sum_{j \neq i}  \frac{1}{|\mathcal{B}_j|} \int_{\mathcal{B}_j}     \frac{x_i-y}{|x_i - y|^3} \otimes (x_i - x_j)    dy    \\
& = \int_{\R^3} \frac{x_i-y}{|x_i - y|^3} \otimes g_n(x_i, y) dy + O(n^{-2/3}) 
\end{align*}
where $g_n(x,y) = \frac{1}{n}\sum_j  \frac{1}{|\mathcal{B}_j|} (x - x_j)  1_{B_j}(y)$ converges for all $x$ weakly * to $g(x,y) = (x-y) \rho(y)$ in $L^\infty_y$ . Thus, $f_n(x) = \int_{\R^3}   \frac{x-y}{|x - y|^3} \otimes g_n(x,y) dy$ converges to  $f(x) = \int_{\R^3}   \frac{x-y}{|x - y|^3} \otimes g(x,y) dy$ and one can check that $f_n$ is equicontinuous over $K$. We conclude as in the case of the Laplacian. 

\mspace
We now turn to the discussion of assumptions \eqref{A1}-\eqref{A2}. We shall see that they are satisfied under the strong assumption  \eqref{minimal_distance} on the minimal distance. The point is  to check \eqref{A2}.  For all $i$ and smooth $\varphi$, we set 
$R_i(x) = \frac{1}{n}\sum_{j \neq i}  \na V(x-x_j) \na \varphi(x_j)$. We shall prove the stronger statement: 
$$ \limsup_n \sum_i \|R_i\|_{L^2(B_i)}^2 \le C \lambda^{1-\eps}  \|\na \varphi\|_{L^\infty}^2 \quad \forall \eps > 0 $$
We have for all $x \in B_i$, 
\begin{equation*}
\begin{aligned}
|R_i(x)| & \le |R_i(x) - R_i(x_i)| + |R_i(x_i)| \le \|\na R_i\|_{L^\infty}|x-x_i| +  |R_i(x_i)| \\
& \le \frac{C}{n} \sum_{j \neq i} \frac{r_n}{|x_i-x_j|^4} \|\na \varphi\|_{L^\infty}+  |R_i(x_i)|  \\
\end{aligned}
\end{equation*} 
By assumption \eqref{minimal_distance}, setting  $y_k = n^{1/3} x_k$, we find
\begin{equation}
 \frac{1}{n} \sum_{j \neq i} \frac{r_n}{|x_i-x_j|^4} \|\na \varphi\|_{L^\infty}  \le   r_n n^{1/3}  \sum_{j \neq i} \frac{1}{|y_i-y_j|^4} \|\na \varphi\|_{L^\infty} \le C \lambda^{1/3} \|\na \varphi\|_{L^\infty}.
\end{equation}  
Hence,  
$$ \sum_i \|a_i\|_{L^2(B_i)}^2 \le C \lambda^{\frac{5}{3}}  \|\na \varphi\|_{L^\infty}^2 + r_n^3 \sum_i |R_i(x_i)|^2.$$
From  a slight variation of \cite[Lemma 2.4]{DGV_MH}, see also \cite{HiWu}, we get for all $q \ge 2$, and  $p$ the conjugate exponent of $q$:
$$ \sum_i |R_i(x_i)|^q \le \frac{C}{\lambda^q} \lambda^{q/p}  \sum_i |\na \varphi(x_i)|^q  =  \frac{C}{\lambda}  \sum_i |\na \varphi(x_i)|^q $$ 
Then, by H\"older inequality: 
\begin{align*}  
r_n^3 \sum_i |R_i(x_i)|^2 & \le  r_n^3 n^{\frac{q-2}{q}} \Big( \sum_i |R_i(x_i)|^q \Big)^{\frac{2}{q}}  \le C r_n^3 n^{\frac{q-2}{q}}  \Big( \frac{1}{\lambda}  \sum_i |\na \varphi(x_i)|^q \Big)^{\frac{2}{q}} \le  C \lambda^{1-\frac{2}{q}}   \|\na \varphi\|_{L^\infty}^2 
\end{align*}
 The result follows.

\mspace
Eventually, it is interesting to understand the meaning of \eqref{A2} when the centers $x_{i,n}$ of the balls are given by a stationary point process  \cite{MR1950431,MR2371524}. More precisely, we will discuss under what conditions on the process we have: for all smooth $\varphi$, 
\begin{equation} \label{A2weak}  \tag{A2$_\flat$}
 \limsup_{n \rightarrow +\infty} \frac{1}{n^2}\sum_i \int_{B_i} \Big| \sum_{j \neq i} \na V(x-x_j) \na \varphi(x_j)\Bigr|^2 dx \  \rightarrow \ 0, \quad \lambda \rightarrow 0. 
\end{equation}
This assumption is slightly weaker than \eqref{A2}:  it only  implies that for  $\varphi \in \dot{H}^1 \cap L^6 \cap C^\infty$,  $\langle R_\lambda , \varphi \rangle = o(\lambda)$,  with $R_\lambda$ the remainder in \eqref{sys:u}. 

\mspace
Let  $\Lambda$ a random point process on a probability space $\Omega$. In particular, for $\omega \in \Omega$, $\Lambda(\omega)$ is a discrete subset of $\R^3$. We assume that the process is stationary, of mean intensity $\lambda_0$, and ergodic. Note that we allow $\lambda_0$ to depend on $\lambda$. Then,  given a small parameter $\eps$, we set $\{x_{1,n}, \dots, x_{n,n}\} = \eps \Lambda \cap \mO$, where the labeling of the centers is arbitrary. Note that  $n$ depends on $\eps$, and is random:  by the ergodic theorem, 
$n \eps^3 \rightarrow \lambda_0 |\mO|$ almost surely   as $\eps \rightarrow 0$. It implies that $r_n \, \eps^{-1} \rightarrow ( \frac{3}{4\pi} \lambda/\lambda_0)^{1/3}$ almost surely   as $\eps \rightarrow 0$. It allows to reformulate condition \eqref{A2weak}: for all smooth $\varphi$, 
$$ \limsup_{\eps \rightarrow 0} \frac{\eps^6}{\lambda_0^2} \sum_i \int_{B_i} \Big| \sum_{j \neq i} \na V(x-x_j) \na \varphi(x_j)\Bigr|^2 dx   \rightarrow 0, \quad \lambda \rightarrow 0. $$
We want to  understand under what conditions  
\begin{equation*} 
 \e  \limsup_{\eps \rightarrow 0}\frac{\eps^6}{\lambda_0^2}  \sum_i \int_{B_i} \Big| \sum_{j \neq i} \na V(x-x_j) \na \varphi(x_j)\Bigr|^2 dx   \rightarrow 0, \quad \lambda \rightarrow 0 
\end{equation*}
which  implies \eqref{A2weak} in probability. First, we want  to reverse the limsup  and the expectation. This will follow from the dominated convergence theorem if we show an $L^\infty$  bound on 
$$  I_{\eps} = \frac{\eps^6}{\lambda_0^2}  \sum_i \int_{B_i} \Big| \sum_{j \neq i} \na V(x-x_j) \na \varphi(x_j)\Bigr|^2 dx $$
that is uniform in $\eps$ and $\omega$ (but not necessarily on $\lambda$). We claim that such bound holds if for instance the process satisfies the condition 
\begin{equation} \label{inf_y_i}
 |y_i - y_j| \ge c (\lambda/\lambda_0)^{1/3}, \quad c > 2 \quad \forall i \neq j
 \end{equation}
which ensures that for $\eps$ small enough, assumption \eqref{A1} holds.  Indeed, using that $\na V$ is a harmonic function, we can write
\begin{align*}
 I_{\eps} & = \frac{\eps^6}{\lambda_0^2} \sum_i \int_{B_i} \Big| \sum_{j \neq i} \dashint_{B_j} \na V(x-x_j) dy \na \varphi(x_j) \Bigr|^2 dx  =  \frac{\eps^6}{\lambda_0^2 r_n^6} \sum_i \int_{B_i} \Big|  \na^2 \Delta^{-1}  \sum_{j \neq i} 1_{B_j}  \na \varphi(x_j)  \Bigr|^2 dx.
 \end{align*}
Note that $|\frac{\eps^6}{r_n^6}| \le \frac{C_0}{\lambda^2}$ for an absolute constant $C_0$. We also remind that $\lambda_0$ possibly depends on $\lambda$. Hence, we find that 
 $$ I_{\eps} \le C_\lambda  \sum_i \int_{B_i} |\na^2 \Delta^{-1} 1_{B_i} \na \varphi(x_i)|^2 + C_\lambda \int_{\R^3} | \na^2 \Delta^{-1}\sum_j 1_{B_j} \na \varphi(x_j)|^2 \le   C'_\lambda $$
 where we use the well-know fact that $\na^2 \Delta^{-1} 1_{B_i}$ is bounded in $L^\infty$ to control the first term,  while we use the continuity of $\na^2 \Delta^{-1}$ over $L^2$ and the fact that the balls $B_j$ are disjoint to control the second term.  
Hence, it remains to understand under what conditions one has
\begin{equation*} 
  \limsup_{\eps \rightarrow 0} \e \, \frac{\eps^6}{\lambda_0^2}  \sum_i \int_{B_i} \Big| \sum_{j \neq i} \na V(x-x_j) \na \varphi(x_j)  \Bigr|^2 dx   \rightarrow 0, \quad \lambda \rightarrow 0. 
\end{equation*}
Let $\tilde B_i = B(x_i, c_0 (\lambda/\lambda_0)^{1/3} \eps)$, with $c_0 = (\frac{3}{4\pi})^{1/3}$, so that almost surely $r_n \sim c_0 (\lambda/\lambda_0)^{1/3} \eps$ as $\eps \rightarrow 0$. Introducing 
$$ \tilde I_\eps = \frac{\eps^6}{\lambda_0^2} \sum_i \int_{\tilde B_i} \Big| \sum_{j \neq i} \na V(x-x_j) \na \varphi(x_j)  \Bigr|^2 dx$$
 we show as in the case of $I_\eps$ that it is bounded uniformly in $\eps$ and $\omega$, and moreover that 
\begin{align*} 
 | I_{\eps} - \tilde I_\eps| & \le C_\lambda \sum_i \int_{B_i \Delta \tilde B_i} |\na^2 \Delta^{-1} 1_{B_i}  \na \varphi(x_i)|^2 + 
 C_\lambda \int_{\cup (B_i \Delta \tilde B_i)} | \na^2 \Delta^{-1}\sum_j 1_{B_j}  \na \varphi(x_j) |^2 \\
 & \le C'_\lambda   \, \Big|1 -  \frac{c_0 \lambda^{1/3} \eps}{\lambda_0^{1/3}r_n}\Big|  \: + \: C_\lambda \|\sum_{i} 1_{B_i \Delta \tilde B_i}\|_{L^2}\|\sum_j 1_{B_j} \na \varphi(x_j)\|_{L^4}^{2} \\
 & \le C''_\lambda   \, \Big(  \big|1 -  \frac{c_0 \lambda^{1/3} \eps}{\lambda_0^{1/3}r_n}\Big|  \: + \:  \Big|1 -  \frac{c_0 \lambda^{1/3} \eps}{\lambda_0^{1/3} r_n}\Big|^{1/2} \Big) \rightarrow 0, \quad \text{as } \eps \rightarrow 0. 
\end{align*}
We used the continuity of $\na^2 \Delta^{-1}$ over $L^4$ in the second  inequality. By the dominated convergence theorem,  the final step is to understand when
\begin{equation*} 
  \limsup_{\eps \rightarrow 0} \e \frac{\eps^6}{\lambda_0^2}  \sum_i \int_{\tilde B_i} \Big| \sum_{j \neq i} \na V(x-x_j) \na \varphi(x_j)\Bigr|^2 dx   \rightarrow 0, \quad \lambda \rightarrow 0.
\end{equation*}
The advantage of $\tilde B_i$ over $B_i$ is that its radius is not random anymore. We write
\begin{equation} 
\begin{aligned}
& \e \frac{\eps^6}{\lambda_0^2} \sum_i \int_{\tilde B_i} \Big| \sum_{j \neq i}\na V(x-x_j) \na \varphi(x_j)\Bigr|^2 dx = \e \frac{\eps^3}{\lambda_0^2} \sum_i \int_{\eps^{-1}\tilde B_i} \Big| \sum_{j \neq i} \na V(y-y_j)  \na \varphi(x_j) \Bigr|^2 dy  \\
= \: &  \e \frac{\eps^3}{\lambda_0^2}  \sum_{i \neq j} \int_{\eps^{-1}\tilde B_i} \Big|\na V(y-y_j)  \na \varphi(x_j)  \Bigr|^2 dy \\ 
+ \: &  \e \frac{\eps^3}{\lambda_0^2}  \sum_{i \neq j \neq j'} \int_{\eps^{-1}\tilde B_i}  \na V(y-y_j)\na \varphi(x_j) \cdot   \na V(y-y_{j'}) \na \varphi(x_{j'}) dy.
 \end{aligned}
\end{equation}
Using the definition of the $k$-point correlation functions $\rho_k$, see \cite[p18]{MR3046995}, we find 
\begin{align*} 
&\:  \e \frac{\eps^3}{\lambda_0^2}  \sum_{i \neq j} \int_{\eps^{-1}\tilde B_i} \Big|\na V(y-y_j)  \na \varphi(x_j)  \Bigr|^2 dy \\ 
 = \: & \: \frac{\eps^3}{\lambda_0^2}   \int_{(\eps^{-1} \mO)^2} \hspace{-0.5cm} dz dz' \rho_2(z,z') \int_{B(z, c_0 (\lambda/\lambda_0)^{1/3})}  \Big|\na V(y-z') \na \varphi(\eps z')\Bigr|^2 dy \\
 \le \: & \: C \frac{\eps^3}{\lambda_0^2}  \int_{(\eps^{-1} \mO)^2} \hspace{-0.5cm} dz dz' \rho_2(z,z')  \int_{B(z, c_0 (\lambda/\lambda_0)^{1/3})}  \frac{1}{|y-z'|^6} dy 
\end{align*}
By stationarity,  $\rho_2(z,z') = \rho_2(0,z'-z)$. After a change of variable, we get
\begin{align*} 
\e \frac{\eps^3}{\lambda_0^2}  \sum_{i \neq j} \int_{\tilde B_i} \Big|\na V(y-y_j)  \na \varphi(x_j)  \Bigr|^2 dy \le   C  \frac{|\mO|}{\lambda_0^2}  \int_{\R^3} \int_{B(0, c_0 (\lambda/\lambda_0)^{1/3})} \frac{\rho_2(0,z')}{|y-z'|^6} dy dz'. 
\end{align*}
Under assumption \eqref{inf_y_i} we find eventually that the r.h.s. goes to zero if and only if
\begin{equation}  \label{cond_rho2}
\frac{\lambda}{\lambda_0^3}   \int_{\R^3} \frac{\rho_2(0,z')}{|z'|^6 + (\lambda/\lambda_0)^2} dz' \rightarrow 0.
\end{equation}
As regards the other term, we get
\begin{align*}
& \e \frac{\eps^3}{\lambda_0^2}  \sum_{i \neq j \neq j'} \int_{\eps^{-1}\tilde B_i}  \na V(y-y_j)\na \varphi(x_j) \cdot   \na V(y-y_{j'}) \na \varphi(x_{j'}) dy \\
= \: & \e \frac{\eps^3}{\lambda_0^2} \int_{(\eps^{-1} \mO)^3} \hspace{-0.6cm} dzdz'dz'' \rho_3(z,z',z'') \int_{B(z, c_0 (\lambda/\lambda_0)^{1/3})} \na V(y-z') \na \varphi(\eps z') \cdot  \na V(y-z'') \na \varphi(\eps z'') dy \\
= \: &  \e \frac{\eps^3}{\lambda_0^2} \int_{(\eps^{-1} \mO)^3} \hspace{-0.6cm} dzdz'dz'' \rho_3(0,z',z'') \int_{B(0, c_0 (\lambda/\lambda_0)^{1/3})} \na V(y-z') \na \varphi(\eps (z + z')) \cdot  \na V(y-z'') \na \varphi(\eps (z +  z'')) dy 
\end{align*}
where we have used stationarity:  $\rho_3(z,z',z'') = \rho_3(0,z'-z, z''-z)$ and a change of variable to obtain the last line. Finally, 
\begin{align*}
& \Big| \e \frac{\eps^3}{\lambda_0^2}  \sum_{i \neq j \neq j'} \int_{\eps^{-1}\tilde B_i}  \na V(y-y_j)\na \varphi(x_j) \cdot   \na V(y-y_{j'}) \na \varphi(x_{j'}) dy \Big|\\
\le \:  & |\mO| \sup_{x \in \R^3}  \int_{B(0, c_0 (\lambda/\lambda_0)^{1/3})}   \Big(\frac{\pa}{\pa z'_i}  \frac{\pa}{\pa z'_k}   \Delta_{z'}^{-1} \Big)  \Big(\frac{\pa}{\pa z''_i}  \frac{\pa}{\pa z''_j}\Delta_{z''}^{-1} \Big) A^\eps_{x,k,j}(z'=y, z''=y) dy 
\end{align*}
where $ A^\eps_{x,k,j}(z', z'')  = 1_{\eps^{-1} \mO}(z') \,  1_{\eps^{-1} \mO}(z'') \, \rho_3(0,z',z'') \, \pa_k \varphi(x + \eps z') \, \pa_j \varphi(x + \eps z'').$

\mspace
Although this last quantity is a bit intricate, we may expect that it goes to zero with $\lambda$ under reasonable assumptions on the process. For instance, if the process is Poisson of intensity $\lambda_0$ (neglecting previous assumptions on the minimal distance), then $\rho_3 = \lambda_0^3$, and  
$|A^\eps_{x,k,j}(z', z'')| \le C \lambda_0^3$, resulting in 
$$ \Big| \e \frac{\eps^3}{\lambda_0^2}  \sum_{i \neq j \neq j'} \int_{\eps^{-1}\tilde B_i}  \na V(y-y_j)\na \varphi(x_j) \cdot   \na V(y-y_{j'}) \na \varphi(x_{j'}) dy \Big| = O(\lambda). $$

\mspace
The condition \eqref{cond_rho2} is more stringent: for the Poisson process of intensity $\lambda_0$, where $\rho_2 = \lambda_0^2$, we find that the quantity at the left-hand side of \eqref{cond_rho2} is  $O(1)$ as $\lambda \rightarrow 0$, but non-vanishing. Hence, our result does not cover this case: to cover it,  we would need a weaker criterion than \eqref{A2}, in the same way as the criterion \eqref{weakseparation} derived in  \cite{MR609184} is weaker than  \eqref{H2} or \eqref{H2sharp}.  Still, \eqref{cond_rho2}  is fulfilled by much more configurations than those satisfying \eqref{minimal_distance} almost surely. Indeed, this latter case corresponds to $\rho_2(0,z') = 0$ for $|z'| \ge c \lambda_0^{-1/3}$, so that the quantity at the left-hand side of \eqref{cond_rho2} is $O(\lambda)$, much stronger than the $o(1)$ asked in \eqref{A2weak}.   

\section*{Acknowledgements}
The author acknowledges the support of the Institut Universitaire de France, and of the SingFlows project, grant ANR-18-CE40-0027 of the French National Research Agency (ANR).

{\footnotesize
 \bibliographystyle{abbrv}
 \bibliography{effective_viscosity_refs}
}

\end{document}